\documentclass{elsart}
\usepackage{amssymb,amsmath}
\usepackage{latexsym}
\journal{Symmetry, Integrability and Geometry: Methods and
Applications}

\begin{document}
\begin{frontmatter}

\title{Solutions for one class of nonlinear fourth-order partial differential equations}

\author[Thailand]{S. Suksern}
\ead{supapornsu@nu.ac.th}

\address[Thailand]{Department of Mathematics, Faculty of Science, Naresuan University,
Phitsanulok, 65000, Thailand}

\begin{abstract}
Some solutions for one class of nonlinear fourth-order partial
differential equations
\[u_{tt}  =
\left( {\kappa u + \gamma u^2 } \right)_{xx}  + \nu uu_{xxxx}  + \mu
u_{xxtt}  + \alpha u_x u_{xxx}  + \beta u_{xx}^2 \] where $\alpha
,\;\beta ,\;\gamma ,\;\mu ,\,\nu $    and $\kappa $ are arbitrary
constants are presented in the paper. This equation may be thought
of as a fourth-order analogue of a generalization of the
Camassa-Holm equation, about which there has been considerable
recent interest. Furthermore, this equation is a Boussinesq-type
equation which arises as a model of vibrations of harmonic
mass-spring chain. The idea of travelling wave solutions and
linearization criteria for fourth-order ordinary differential
equations by point transformations are applied to this problem.

\end{abstract}

\end{frontmatter}

\section{Introduction}
\label{Intro}

Almost all important governing equations in physics take the form of
nonlinear differential equations, and, in general, are very
difficult to solve explicitly.  While solving problems related to
nonlinear ordinary differential equations it is often expedient to
simplify equations by a suitable change of variables.

 Many methods of solving differential equations use a
change of variables that transforms a given differential equation
into another equation with known properties. Since the class of
linear equations is considered to be the simplest class of
equations, there is the problem of transforming a given differential
equation into a linear equation. This problem, which is called a
linearization problem. The reduction of an ordinary differential
equation to a linear ordinary differential equation besides
simplification allows constructing an exact solution of the original
equation.

One of the most interested nonlinear problems but also difficultly
in solving is the problem of nonlinear fourth-order partial
differential equations \cite{bk:ClarksonMansfieldPriestley[1997]}
\begin{equation}\label{sc:12}
    u_{tt}  = \left( {\kappa u + \gamma u^2 }
\right)_{xx}  + \nu uu_{xxxx}  + \mu u_{xxtt}  + \alpha u_x u_{xxx}
+ \beta u_{xx}^2
\end{equation}
where $\alpha ,\;\beta ,\;\gamma ,\;\mu ,\,\nu $
  and  $\kappa $
are arbitrary constants. The main difficulty in solving this problem
comes from the terms of nonlinear partial differential equations and
the large number of order. Because of this difficulty, there are
only a few attempts to solve this problem.

In 2008, Suksern, Meleshko and Ibragimov
\cite{bk:IbragimovMeleshkoSuksern[2008],bk:IbragimovMeleshkoSuksern[2007]}
found the explicit form of the criteria for linearization of
fourth-order ordinary differential equations by point
transformations. Moreover, the procedure for the construction of
the linearizing transformation are presented. By the virtue of
\cite{bk:IbragimovMeleshkoSuksern[2008],bk:IbragimovMeleshkoSuksern[2007]}
to bring about the idea for solving the problem of nonlinear
fourth-order partial differential equations (\ref{sc:12}).

The way of solving the problem is organized as follows. Firstly,
reducing the nonlinear fourth-order partial differential equations
to the nonlinear fourth-order ordinary differential equations by
substituting the form of travelling wave solutions. Secondly,
reducing the nonlinear fourth-order ordinary differential equations
to the linear fourth-order ordinary differential equations by
applying the criteria for linearization in
\cite{bk:IbragimovMeleshkoSuksern[2008],bk:IbragimovMeleshkoSuksern[2007]}.
Finally, finding the exact solutions of linear equations and then
substituting back to the exact solutions of the original problem.

\section{Linearization criteria for fourth-order ordinary differential equations by point transformations}
\label{point}

The important tools for this research is the linearization
criteria for fourth-order ordinary differential equations by point
transformations. From
\cite{bk:IbragimovMeleshkoSuksern[2008],bk:IbragimovMeleshkoSuksern[2007]}
we have the following theorems.

\begin{thm}
Any fourth-order ordinary differential equation
$$ y^{\left( 4
\right)}  = f\left( {x,y,y',y'',y'''} \right),$$
can be reduced by a
point transformation
\begin{equation}\label{sc:02}
            t = \varphi \left( {x,y} \right),\quad u = \psi \left( {x,y}
            \right),
        \end{equation}
to the linear equation
\begin{equation}
 \label{sc:01}
u^{\left( 4 \right)}  + \alpha \left( t \right)u' + \beta \left( t
\right)u = 0,
 \end{equation}
  where $t$ and $u$ are  the independent and dependent variables,
  respectively,
if it belongs to the class of equations
 \begin{equation}
 \label{sc:03}
 \begin{array}{ll}
  y^{(4)}  & + (A_{1}y' + A_{0})y''' + B_{0}y''^2 + (C_{2}y'^2 + C_{1}y' + C_{0})y''
 \\ &+ D_{4}y'^4 + D_{3}y'^3 + D_{2}y'^2 + D_{1}y' + D_{0} = 0,
\end{array}
 \end{equation}
 or
 \begin{equation}
 \label{sc:04}
\begin{array}{ll}
  y^{( 4)}   &+  \frac{1}{y'  +  r}
(  - 10 y''  +  F_2 y'^2  +  F_1 y'  +  F_{0 } )y'''
   \\&+ \frac{1}
{{( y'  +  r)^2 }}\left[
  15y''^3   +  ( H_2 y'^2  +  H_1 y' +  H_0 )y''^2  \right.
  \\& +  (J_4 y'^4  +  J_3 y'^3  +  J_2 y'^2  +  J_1 y'  + J_{0 } )y'' 
   \\&+  K_7 y'^7  +  K_6 y'^6  +  K_5 y'^5  +  K_4 y'^4   \\&
\left.   +  K_3 y'^3  +  K_2 y'^2  +  K_1 y' +  K_0
  \right] = 0, 
\end{array}
\end{equation}
 where $A_j = A_j(x,y)$, $B_j = B_j(x,y)$, $C_j = C_j(x,y)$, $ D_j =
 D_j(x,y)$, $r=r(x,y), F_j = F_j(x,y)$, $H_j = H_j(x,y)$, $J_j = J_j(x,y)$ and $ K_j =
 K_j(x,y)$  are arbitrary functions of $x, y$.
\end{thm}

Since this research deals with the first class, let us emphasize to
the first class for other theorems that we need to use.

\begin{thm}
 \label{eq:09th}
  Equation (\ref{sc:03}) is linearizable if and only if its coefficients obey the following
conditions:
\begin{align}
& A_{0y} - A_{1x} = 0,\label{eq:45}\\[1.5ex]
& 4B_{0} - 3A_{1} = 0, \label{eq:46}\\[1.5ex]
& 12A_{1y} + 3A_{1}^2 - 8C_{2} = 0, \label{eq:47}\\[1.5ex]
& 12A_{1x} + 3A_{0}A_{1} - 4C_{1} = 0, \label{eq:48}\\[1.5ex]
& 32C_{0y} + 12A_{0x}A_{1} - 16C_{1x} + 3A_{0}^2A_{1} - 4A_{0}C_{1}
= 0, \label{eq:49}\\[1.5ex]
& 4C_{2y} + A_{1}C_{2} - 24D_{4} = 0, \label{eq:50}\\[1.5ex]
& 4C_{1y} + A_{1}C_{1} - 12D_{3} = 0, \label{eq:51}\\[1.5ex]
& 16C_{1x} - 12A_{0x}A_{1} - 3A_{0}^2A_{1} + 4A_{0}C_{1} +
8A_{1}C_{0} - 32D_{2}=0, \label{eq:52}\\[1.5ex]
\nonumber
 & 192D_{2x} + 36A_{0x}A_{0}A_{1} - 48A_{0x}C_{1} -
48C_{0x}A_{1} - 288D_{1y} + 9A_{0}^3A_{1} \nonumber \\ & -
12A_{0}^2C_{1} - 36A_{0}A_{1}C_{0} + 48A_{0}D_{2} + 32C_{0}C_{1} =
0, \label{eq:53}
\end{align}
\begin{align}
\nonumber
 & 384D_{1xy} - \Big[3((3A_{0}A_{1} - 4C_{1})A_{0}^2 +
16(2A_{1}D_{1} + C_{0}C_{1}) \nonumber \\
& - 16(A_{1}C_{0} - D_{2})A_{0})A_{0} - 32(4(C_{1}D_{1} -
2C_{2}D_{0} + C_{0}D_{2}) \nonumber \\
& + (3A_{1}D_{0} - C_{0}^2)A_{1}) - 96D_{1y}A_{0} + 384D_{0y}A_{1}
+ 1536D_{0yy} \nonumber \\
& - 16(3A_{0}A_{1} - 4C_{1})C_{0x} + 12((3A_{0}A_{1} -
4C_{1})A_{0} \nonumber \\
& - 4(A_{1}C_{0}  - 4D_{2}))A_{0x}\Big] = 0. \label{eq:54}
\end{align}
\end{thm}

\begin{thm}\label{eq:10th}
Provided that the conditions (\ref{eq:45})-(\ref{eq:54}) are
satisfied, the linearizing transformation (\ref{sc:02}) is defined
by a fourth-order ordinary differential equation for the function
$\varphi (x),$ namely by the Riccati equation
 \begin{equation}
 \label{eq:55}
 40\frac{d \chi}{dx} - 20\chi^2 = 8C_0 -3A_0^2 -12A_{0x},
 \end{equation}
 for
 \begin{equation}
 \label{eq:56}
 \chi = \frac{\varphi_{xx}}{\varphi_x}\,,
 \end{equation}
and by the following integrable system of partial differential
equations for the function  $\psi (x,y)$
\begin{equation}\label{eq:57}
    4\psi_{yy} =  \psi_{y}A_{1},
\ \ \
    4\psi_{xy} =  \psi_{y}(A_{0} + 6\chi),
\end{equation}
and
\begin{align}
 \nonumber
& 1600\psi_{xxxx} =  9600\psi_{xxx}\chi + 160\psi_{xx}( - 12A_{0x} -
3A_{0}^2 - 90\chi^2 + 8C_{0}) \nonumber \\ & +
40\psi_{x}(12A_{0x}A_{0} + 72A_{0x}\chi - 16C_{0x} + 3A_{0}^3 +
18A_{0}^2\chi  - 12A_{0}C_{0} \nonumber \\ & + 120\chi^3 - 48\chi
C_{0} + 24D_{1} - 8\Omega)  + \psi(144A_{0x}^2 + 72A_{0x}A_{0}^2 -
352A_{0x}C_{0} \nonumber \\ &- 160C_{0xx} - 80C_{0x}A_{0}
 - 1600D_{0y} + 640D_{1x} - 80\Omega_{x} + 9A_{0}^4
- 88A_{0}^2C_{0} \nonumber \\ & + 160A_{0}D_{1}  + 30A_{0}\Omega -
400A_{1}D_{0} + 300\chi\Omega + 144C_{0}^2) + 1600\psi_{y}D_{0},
\label{eq:59}
\end{align}
where $\chi$ is given by equation (\ref{eq:56}) and $\Omega$ is the
 following expression
\begin{equation}
 \label{eq:60}
 \Omega = A_{0}^3 - 4A_{0}C_{0} + 8D_{1} - 8C_{0x}
+ 6A_{0x}A_{0} + 4A_{0xx}.
\end{equation}
Finally, the coefficients $\alpha $ and $\beta$ of the resulting
linear equation (\ref{sc:01}) are
\begin{equation}\label{eq:61}
    \alpha =  \frac{\Omega}{8\varphi_{x}^3} ~,
\end{equation}
\begin{align}
\nonumber
 & \beta =  (1600\varphi_{x}^4)^{-1}( - 144A_{0x}^2
- 72A_{0x}A_{0}^2 + 352A_{0x}C_{0} + 160C_{0xx} + 80C_{0x}A_{0}
\nonumber \\ & + 1600D_{0y} - 640D_{1x} + 80\Omega_{x} - 9A_{0}^4 +
88A_{0}^2C_{0} - 160A_{0}D_{1} - 30A_{0}\Omega \nonumber \\ &+
400A_{1}D_{0} - 300\chi\Omega - 144C_{0}^2).\label{eq:62}
\end{align}
\end{thm}
\section{Method and Result}
\label{result}

Let us consider the nonlinear fourth-order partial differential
equation (Clarkson and Priestley, 1999)
\begin{equation}\label{pde:01}
u_{tt}  = \left( {\kappa u + \gamma u^2 } \right)_{xx}  + \nu
uu_{xxxx}  + \mu u_{xxtt}  + \alpha u_x u_{xxx}  + \beta u_{xx} ^2 ,
\end{equation}
where $\alpha, \beta, \gamma, \mu, \nu$ and $\kappa$ are arbitrary
constants.

Of particular interest among solutions of equation (\ref{pde:01})
are traveling wave solutions:
\begin{equation*}
    u(x,t)=H(x-Dt),
\end{equation*}
where $D$ is a constant phase velocity and the argument $x-Dt$ is a
phase of the wave. Substituting the representation of a solution
into equation (\ref{pde:01}), one finds
\begin{equation}\label{pde:02}
    (\nu H+\mu D^{2})H^{(4)} + \alpha H'H''' + \beta H''^{2} + (2\gamma
    H+\kappa-D^{2})H'' + 2\gamma H'^{2}=0.
\end{equation}
This is an equation of the form (\ref{sc:03}) with coefficients
\begin{equation*}
    A_{1} = \frac{\alpha}{\nu H+\mu D^{2}}, \,\,\, A_{0}=0, \,\,\, B_{0} = \frac{\beta}{\nu H+\mu
    D^{2}},
\end{equation*}
\begin{equation*}
    C_{2}=C_{1}=0, \,\,\, C_{0}=\frac{2\gamma H + \kappa - D^{2}}{\nu H+\mu D^{2}},
\end{equation*}
\begin{equation*}
    D_{4}=D_{3}=0, \,\,\, D_{2}=\frac{2\gamma}{\nu H+\mu
    D^{2}}, \,\,\, D_{1}=D_{0}=0.
\end{equation*}
Substituting these coefficients into the linearization conditions
(\ref{eq:45})-(\ref{eq:54}), one obtains the following results.
\subsection{Case 1 : $\nu = 0$}

If $\nu = 0$, then equation (\ref{pde:02}) is linearizable if and
only if
\begin{equation*}
    \alpha=\beta=\gamma=0.
\end{equation*}
This means, with theses conditions  the original equation
(\ref{pde:02}) becomes the linear equation
\[
    (\mu D^{2})H^{(4)} +(\kappa-D^{2})H'' =0.
\]
The solution of this equation is
\[
H\left( {x - Dt} \right) = C_1 \sin \sqrt {\frac{{\kappa  - D^2 }}
{{\mu D^2 }}} \left( {x - Dt} \right) + C_2 \cos \sqrt
{\frac{{\kappa  - D^2 }} {{\mu D^2 }}} \left( {x - Dt} \right),
\]
where $C_1$ and $C_2$ are arbitrary constants. Hence,
\[
u\left( {x , t} \right) = C_1 \sin \sqrt {\frac{{\kappa  - D^2 }}
{{\mu D^2 }}} \left( {x - Dt} \right) + C_2 \cos \sqrt
{\frac{{\kappa  - D^2 }} {{\mu D^2 }}} \left( {x - Dt} \right).
\]
\subsection{Case 2 : $\nu \neq 0$}
\textbf{Case 2.1 : $\gamma = 0$ }

\begin{itemize}
    \item If $\nu \neq 0$, $\gamma = 0$ and $\beta = 0$, then equation (\ref{pde:02}) is
linearizable if and only if
\begin{equation*}
    \alpha=0, \,\,\, \kappa=D^2.
\end{equation*}
Next, to find the linear form of equation (\ref{pde:02}) and it's solutions.\\

\emph{Note that} : For applying the theorems to our problem, here
$x=x-Dt$ and $y(x)=H(x-Dt)$. \\

\noindent From (\ref{eq:55}) one has
$$
    2\frac{{d\chi }} {{dx}} - \chi ^2  = 0.
$$
Let us take its simplest solution  $\chi=0$ .  Then invoking
(\ref{eq:56}), we let
$$
   \varphi  = x.
$$
Now the equations (\ref{eq:57})-(\ref{eq:59}) are written as
$$
    \psi _{yy} =0, \quad \psi_{xy}  = 0,
$$
and yield
$$
    \psi _y  = K_0 ,\quad \quad K_0 = const.
$$
Hence,
$$
   \psi  = K_0 y + K_1\left( x \right).
$$
Since one can use any particular solution, we set  $ K_0 =
1,\;\;K_1\left( x \right) = 0$  and take
$$
    \psi  = y.
$$
Noting that (\ref{eq:60}) yields $\Omega = 0 $, one can readily
verify that the function  $ \psi  = y $ solves equation
(\ref{eq:59}) as well. Hence, one obtains the following
transformations
\begin{equation}\label{sc:06}
    \tilde{t} = x,\quad \tilde{u} = y.
\end{equation}
Since  $ \Omega  = 0 $,  equations (\ref{eq:61}) and (\ref{eq:62})
give
$$
   \tilde{\alpha}  = 0,\quad  \tilde{\beta}  = 0.
$$
Hence, the equation (\ref{pde:02}) is mapped by the
transformations (\ref{sc:06}) to the linear equation
$$
\tilde{u}^{\left( 4 \right)}  = 0.
$$
The solution of this linear equation is
$$
\tilde{u}=C_0 + C_1 \tilde{t} + C_2 \tilde{t}^2 +C_3 \tilde{t}^3,
$$
where $C_0, C_1, C_2$ and $C_3$ are arbitrary constants. That is
$$
H(x-Dt)=C_0 + C_1 (x-Dt) + C_2 (x-Dt)^2 +C_3 (x-Dt)^3.
$$
Hence,
$$
u(x,t)=C_0 + C_1 (x-Dt) + C_2 (x-Dt)^2 +C_3 (x-Dt)^3.
$$\\
    \item If $\nu \neq 0$, $\gamma = 0$ and $\beta = 3\nu$, then equation (\ref{pde:02}) is
linearizable if and only if
\begin{equation*}
    \alpha=4\nu, \,\,\,  \kappa=D^2 .
\end{equation*}
\end{itemize}
Next, to find the linear form of equation (\ref{pde:02}) and it's
solutions. Considering (\ref{eq:55}) one has
$$
    2\frac{{d\chi }} {{dx}} - \chi ^2  = 0.
$$
Let us take its simplest solution  $\chi=0$ .  Then invoking
(\ref{eq:56}), we let
$$
   \varphi  = x.
$$
Now the equations (\ref{eq:57})-(\ref{eq:59}) are written as
$$
    \psi _{yy} =\frac{{\nu \psi _y }} {{\left( {\nu y + D^2 \mu } \right)}}
, \quad \psi_{xy}  = 0,
$$
and yield
\[
\psi _y  = K_0 \left( x \right)\left( {y + \frac{{D^2 \mu }} {\nu }}
\right).
\]
Since $\psi_{xy}  = 0$, that means $K_0 (x) =K_0=const$. One arrives
at
\[
\psi  = K_0 \left( {\frac{{y^2 }} {2} + \frac{{D^2 \mu }} {\nu }y}
\right) + K_1 (x).
\]
Since one can use any particular solution, we set  $ K_0 =
1,\;\;K_1\left( x \right) = 0$  and take
$$
    \psi  = {\frac{{y^2 }} {2} + \frac{{D^2 \mu }} {\nu }y}.
$$
Observing that (\ref{eq:60}) yields $\Omega = 0 $, one can readily
verify that the function  $ \psi  = {\frac{{y^2 }} {2} + \frac{{D^2
\mu }} {\nu }y} $ solves equation (\ref{eq:59}) as well. Hence, one
obtains the following transformations
\begin{equation}\label{sc:08}
    \tilde{t} = x,\quad \tilde{u} = {\frac{{y^2 }} {2} + \frac{{D^2 \mu }} {\nu }y}.
\end{equation}
Since  $ \Omega  = 0 $,  equations (\ref{eq:61}) and (\ref{eq:62})
give
$$
   \tilde{\alpha}  = 0,\quad  \tilde{\beta}  = 0.
$$
Hence, the equation (\ref{pde:02}) is mapped by the
transformations (\ref{sc:08}) to the linear equation
$$
\tilde{u}^{\left( 4 \right)}  = 0.
$$
The solution of this linear equation is
$$
\tilde{u}=C_0 + C_1 \tilde{t} + C_2 \tilde{t}^2 +C_3 \tilde{t}^3,
$$
where $C_0, C_1, C_2$ and $C_3$ are arbitrary constants. So that we
obtain the implicit solution in the form
$$
{\frac{{H^2 }} {2} + \frac{{D^2 \mu }} {\nu }H}=C_0 + C_1 (x-Dt) +
C_2 (x-Dt)^2 +C_3 (x-Dt)^3.
$$
Hence,
$$
{\frac{{u^2 }} {2} + \frac{{D^2 \mu }} {\nu }u}=C_0 + C_1 (x-Dt) +
C_2 (x-Dt)^2 +C_3 (x-Dt)^3.
$$\\

\textbf{Case 2.2 : $\gamma \neq 0$} \\
Equation (\ref{pde:02}) is linearizable if and only if
\begin{equation*}
    \alpha=4\nu, \,\,\, \beta=3\nu, \,\,\, \kappa=\frac{(2\gamma \mu + \nu)D^{2}}{\nu}.
\end{equation*}
Because of (\ref{eq:55}) one obtains
$$
    \frac{{d\chi }} {{dx}} -\frac{1}{2} \chi ^2  = \frac{2\gamma}{5\nu}.
$$
Solving this equation, one gets
\[
\chi  = 2\sqrt {\frac{\gamma } {{5\nu }}} \tan \left( {\sqrt
{\frac{\gamma } {{5\nu }}} \left( {x + C} \right)} \right),
\]
where $C$ is an arbitrary constant. Since one can use any
particular solution, we set $ C = 0$, so that
\[
\chi  = 2\sqrt {\frac{\gamma } {{5\nu }}} \tan \left( {\sqrt
{\frac{\gamma } {{5\nu }}}  {x }} \right).
\]
Then invoking (\ref{eq:56}),
\[
\frac{{\varphi _{xx} }} {{\varphi _x }} = 2\sqrt {\frac{\gamma }
{{5\nu }}} \tan \left( {\sqrt {\frac{\gamma } {{5\nu }}} x}
\right).
\]
Thus,
\[
\varphi _x  = K_0 \sec ^2 \left( {\sqrt {\frac{\gamma } {{5\nu }}}
x} \right),
\]
where $K_0$ is an arbitrary constant. By solving this equation,
one obtains
\[
\varphi  = \sqrt {\frac{{5\nu }} {\gamma }} K_0 \tan \left( {\sqrt
{\frac{\gamma } {{5\nu }}} x} \right) + K_1,
\]
where $K_1$ is an arbitrary constant. One can choose $ K_0  =
\sqrt {\frac{\gamma } {{5\nu }}}$ and $K_1 =0$. One arrives at
\[
\varphi  = \tan \left( {\sqrt {\frac{\gamma } {{5\nu }}} x}
\right).
\]
Now the equations (\ref{eq:57})-(\ref{eq:59}) are written as
\begin{equation}\label{sc:09}
\psi _{xy}  = 3\sqrt {\frac{\gamma } {{5\nu }}} \tan \left( {\sqrt
{\frac{\gamma } {{5\nu }}} x} \right)\psi _y,
\end{equation}
and
\begin{equation}\label{sc:10}
\psi _{yy}  = \frac{{\nu \psi _y }} {{D^2 \mu  + \nu y}}.
\end{equation}
Equation (\ref{sc:09}) and (\ref{sc:10}) give
\[
\psi  = \sec ^3 \left( {\sqrt {\frac{\gamma } {{5\nu }}} x}
\right)\left[ {\frac{{y^2 }} {2} + \frac{{D^2 \mu }} {\nu }y}
\right].
\]
One can readily verify that the function  $ \psi$ solves equation
(\ref{eq:59}) as well. Hence, one obtains the following
transformations
\begin{equation}\label{sc:11}
    \tilde{t} = \tan \left( {\sqrt {\frac{\gamma } {{5\nu }}} x}
\right),\quad \tilde{u} = \sec ^3 \left( {\sqrt {\frac{\gamma }
{{5\nu }}} x} \right)\left[ {\frac{{y^2 }} {2} + \frac{{D^2 \mu }}
{\nu }y} \right].
\end{equation}
Equations (\ref{eq:61}) and (\ref{eq:62}) give
$$
   \tilde{\alpha}  = 0,\quad  \tilde{\beta}  =  - 9\cos ^8 \sqrt {\frac{\gamma }
{{5\nu }}} x .
$$
Hence, the equation (\ref{pde:02}) is mapped by the
transformations (\ref{sc:11}) to the linear equation
\begin{equation}\label{sc:13}
    \tilde{u}^{\left( 4 \right)}  -9\cos ^8 \sqrt {\frac{\gamma } {{5\nu
}}} x \tilde{u}= 0.
\end{equation}
Since this is an unsolved linear equation,  so that the result of
this case we obtained only at the linear form of nonlinear equation
(\ref{pde:02}).

\section{Discussion and Conclusion}
\label{con}

In the present work, we found some following solutions for nonlinear
fourth-order partial differential equations (\ref{pde:01}).
\begin{itemize}
    \item If $\nu = \alpha=\beta=\gamma=0$, then the solution of
    (\ref{pde:01}) is
        \[
            u\left( {x , t} \right) = C_1 \sin \sqrt {\frac{{\kappa  - D^2 }}
            {{\mu D^2 }}} \left( {x - Dt} \right) + C_2 \cos \sqrt
            {\frac{{\kappa  - D^2 }} {{\mu D^2 }}} \left( {x - Dt} \right).
        \]
    \item  If $\nu \neq 0$, $\alpha=\beta=\gamma=0$ and $\kappa=D^2$, then the solution of
    (\ref{pde:01}) is
        \[
            u(x,t)=C_0 + C_1 (x-Dt) + C_2 (x-Dt)^2 +C_3 (x-Dt)^3.
        \]
    \item If $\nu \neq 0$, $\alpha=4\nu$, $\beta = 3\nu, \gamma = 0$ and
    $\kappa=D^2$ , then the solution of
    (\ref{pde:01}) is
        \[
            {\frac{{u^2(x,t) }} {2} + \frac{{D^2 \mu }} {\nu }u(x,t)}=C_0 + C_1 (x-Dt) +
            C_2 (x-Dt)^2 +C_3 (x-Dt)^3.
        \]
    \item If $\gamma \neq 0, \alpha=4\nu,  \beta=3\nu$ and $\kappa=\frac{(2\gamma \mu +
    \nu)D^{2}}{\nu}$,  then the linear form of (\ref{pde:01}) is (\ref{sc:13}).
\end{itemize}
An interesting aspect of the results in this paper is that the class
of exact solutions of the original nonlinear problems, which can not
find by the classical methods.

\section{Acknowledgements}
\label{Acknowledgements}

This work was financially supported by Faculty of Science, Naresuan
University. The author wishes to express thanks to Prof.Dr. Sergey
V. Meleshko, Suranaree University of Technology for his guidance
during the work.

\bibliographystyle{elsart-num}
\bibliography{ref}

\end{document}